\documentclass[12pt]{article}
\usepackage{curves}
\usepackage{amsmath, amssymb, lamsarrow}
\usepackage{amssymb}
\def\mineappendix{
        \setcounter{section}{1}
        \setcounter{subsection}{0}
        \def\thesection{\Alph{section}}
        \def\sectionap{\@startsection  {section}{1}{\z@}
                        {-3.5ex plus-1ex minus-.2ex} {0ex plus.2ex}
                        {\reset@font\Large\bf  Appendix:  \, }
                        }
        }
\makeatother
\def\Proclaim #1. #2\par{\bigbreak\noindent{\sc#1.\enspace}{\it#2}\par}

\newcommand{\gwii}[1]{\left< \hspace{-2pt} \left< \, #1 \,
        \right>  \hspace{-2pt} \right>_{0}}

\newcommand{\gwiione}[1]{\left< \hspace{-2pt} \left< \, #1 \,
        \right> \hspace{-2pt} \right>_{1}}

\newcommand{\gwiitwo}[1]{\left< \hspace{-2pt} \left< \, #1 \,
        \right> \hspace{-2pt} \right>_{2}}

\newcommand{\gwih}[2]{ \left< \, #2 \, \right>_{#1}}
\newcommand{\gwiig}[1]{\left< \hspace{-2pt} \left< \, #1 \,
    \right> \hspace{-2pt} \right>_{g}}
\newcommand{\gwiih}[2]{\left< \hspace{-2pt} \left< \, #2 \,
    \right> \hspace{-2pt} \right>_{#1}}

\newcommand{\gwiigE}[1]{\left< \hspace{-2pt} \left< \, #1 \,
    \right> \hspace{-2pt} \right>_{g, \mathbb{E}}}
\newcommand{\gwiihE}[2]{\left< \hspace{-2pt} \left< \, #2 \,
    \right> \hspace{-2pt} \right>_{#1, \mathbb{E}}}

\newcommand{\grava}[1]{\tau_{#1}(\gamma_{\alpha})}
\newcommand{\gravua}[1]{\tau_{#1}(\gamma^{\alpha})}
\newcommand{\gravb}[1]{\tau_{#1}(\gamma_{\beta})}
\newcommand{\gravub}[1]{\tau_{#1}(\gamma^{\beta})}
\newcommand{\gravm}[1]{\tau_{#1}(\gamma_{\mu})}
\newcommand{\gravum}[1]{\tau_{#1}(\gamma^{\mu})}

\newcommand{\gravun}[1]{\tau_{#1}(\gamma^{\nu})}

\newcommand{\ga}{\gamma_{\alpha}}
\newcommand{\gua}{\gamma^{\alpha}}
\newcommand{\gb}{\gamma_{\beta}}
\newcommand{\gub}{\gamma^{\beta}}
\newcommand{\gm}{\gamma_{\mu}}
\newcommand{\gum}{\gamma^{\mu}}
\newcommand{\gn}{\gamma_{\nu}}
\newcommand{\gun}{\gamma^{\nu}}

\newcommand{\vs}{{\cal S}}
\newcommand{\vx}{{\cal X}}

\newcommand{\vw}{{\cal W}}
\newcommand{\vv}{{\cal V}}
\newcommand{\vu}{{\cal U}}

\newcommand{\qp}{\circ}

\newtheorem{lem}{Lemma}[section]

\newtheorem{thm}[lem]{Theorem}
\newtheorem{pro}[lem]{Proposition}

\title{On Certain Vanishing Identities For \\ Gromov-Witten Invariants}
\author{Xiaobo Liu \thanks{Research was partially supported by NSF grant DMS-0505835}}
\date{}

\begin{document}
\maketitle




 Let
$V$ be a compact symplectic manifold and $\{\ga \mid \alpha=1, \cdots N \}$ be a basis
for $H^{*}(V; \mathbb{C})$. We choose $\gamma_{1}$ to be the identity of cohomology ring of $V$. Let
$\gua = \eta^{\alpha \beta} \gb$ with $(\eta^{\alpha \beta})$ representing the inverse
matrix of the Poinc\'{a}re intersection pairing. As a convention, repeated Greek letter
indices are summed
up over their entire ranges from 1 to $N$.
Recall that the {\it big phase space} for Gromov-Witten invariants is
$\prod_{n=0}^{\infty} H^{*}(V; \mathbb{C})$ with standard basis
$\{ \grava{n} \mid \alpha=1, \ldots, N, \,\,\,\,\, n \geq 0\}$.
Let $t_{n}^{\alpha}$ be the coordinates on the big phase space
with respect to the standard basis. The genus-$g$ generating function $F_{g}$ is
a formal power series of $t=(t_{n}^{\alpha})$ with coefficients given by
genus-g Gromov-Witten invariants. Derivatives of
$F_{g}$ with respect to $t_{n_{1}}^{\alpha_{1}}, \ldots, t_{n_{k}}^{\alpha_{k}}$ are denoted by
$\gwiig{\tau_{n_{1}}(\alpha_{1}) \cdots \tau_{n_{k}}(\alpha_{k})}$.

The following conjecture was proposed by K. Liu and H. Xu in \cite{LkX}:

\vspace{6pt}
\noindent {\bf Conjecture A}: {\it For $m > 2g$, genus-$g$
Gromov-Witten invariants of $V$ satisfy the following identity
\[ \sum_{j=0}^{m} (-1)^{j} \gwiig{\grava{j} \,\, \gravua{m-j}} = 0. \]
}

\noindent Note that when $m$ is odd, this conjecture is trivial dual
to the symmetry of the indices. One can see this by simply replacing $j$ by $m-j$
on the left hand side. So we only need to consider the case
when $m$ is even.

In this paper, we prove the following
\begin{thm} \label{thm:A}
Conjecture A is true if $m \geq 3g + \delta_{g, 0}$ for all genera.
\end{thm}
In particular, this theorem implies that Conjecture A is true for
$g \leq 2$.

Another conjecture proposed in \cite{LkX} is the following:

\vspace{6pt}

\noindent
 {\bf Conjecture B}: {\it For $k \geq g$,
 \[ - \sum_{n, \alpha} \tilde{t}_{n}^{\alpha} \gwiig{\grava{n+2k-1}}
    + \frac{1}{2} \sum_{h=0}^{g} \sum_{j=0}^{2k-2} (-1)^{j}
            \gwiih{h}{\grava{j}} \gwiih{g-h}{\gravua{2k-2-j}} =
            \frac{1}{2} \delta_{g, 0} \delta_{k,0} \eta_{\alpha \beta}
                        t_{0}^{\alpha} t_{0}^{\beta} \]
 where $\tilde{t}_{n}^{\alpha} = t_{n}^{\alpha} -\delta_{\alpha, 1} \delta_{n, 1} $.
}

\noindent

In this paper, we also prove the following
\begin{thm} \label{thm:B}
(a) For $m = 2k >2g$ , Conjecture A  for genus $g$ is equivalent to conjecture B for
genus $g+1$.

(b) Conjecture B is true if $2k \geq 3g-1$ for all genera.

(c) Conjecture B is true for genus $g \leq 3$.
\end{thm}
Part (a) of this theorem follows from a result of Faber and Pandharipande \cite{FP}, which will
be discussed in more detail in Section~\ref{sec:A=B}.

The following conjecture which is more general than conjecture B was also proposed in \cite{LkX}:

\vspace{6pt}
\noindent
{\bf Conjecture C}: For all $x_{i}, y_{i} \in H^{*}(V; \mathbb{C})$ and $m \geq 2g-3+a +b$,
\[ \sum_{j \in \mathbb{Z}}
     \sum_{h=0}^{g} (-1)^{j} \gwiih{h}{\grava{j} \prod_{i=1}^{a} \tau_{p_{i}}(x_{i})}
    \gwiih{g-h}{\gravua{m-j} \prod_{i=1}^{b} \tau_{q_{i}}(y_{i}) } = 0. \]
Here $j$ is allowed to be an arbitrary integer. To interpret this conjecture correctly,
one has to use the convention that $\gwih{0,0}{\tau_{-2}(\gamma_{1})} = 1$
and
\[ \gwih{0,0}{\grava{m} \gravb{-1-m}} = (-1)^{{\rm max}(m, -1-m)} \eta_{\alpha \beta},
    \hspace{20pt} m \in \mathbb{Z}. \]
Conjecture B is
the special case of Conjecture C for $a=b=0$.
Note that we will not use this convention in the rest of this paper.
We will reformulate these conjectures
in Theorem~\ref{thm:reform} in order to see their
relations with topological recursion relations. We will also
prove the following vanishing identity for all genera:
\begin{thm} \label{thm:allg}
Assume $g \geq h \geq 0$, $a, b \geq 0$, $m \geq 3g-3+a+b$. Moreover we also assume
 $a \geq 2$ if $h=0$ and $b \geq 2$ if $h=g$. Then for all $p_{i}, q_{j} \geq 0$
 and $x_{i}, y_{j} \in H^{*}(V; \mathbb{C})$, we have
\begin{equation} \label{eqn:allg1}
     \sum_{j=0}^{m}
      (-1)^{j} \gwiih{h}{\grava{j} \prod_{i=1}^{a} \tau_{p_{i}}(x_{i})}
    \gwiih{g-h}{\gravua{m-j} \prod_{i=1}^{b} \tau_{q_{i}}(y_{i}) } = 0.
\end{equation}
In the case that $h=0$, $g \geq 1$ and $a=1$, we also have the following identity
\begin{equation} \label{eqn:allg2}
     \sum_{j=0}^{m}
      (-1)^{j} \gwii{\grava{j} \tau_{p}(x)}
    \gwiig{\gravua{m-j} \prod_{i=1}^{b} \tau_{q_{i}}(y_{i}) }
    = \gwiig{\tau_{p+m+1}(x) \prod_{i=1}^{b} \tau_{q_{i}}(y_{i}) }
\end{equation}
for $m \geq 3g-2+b$.
\end{thm}
Note that the vanishing identity in this theorem is stronger than corresponding cases for
Conjecture C since there is no summation over genus.

This paper is organized as follows: In Section~\ref{sec:A=B},
we will discuss a result in \cite{FP} which implies part (a) of Theorem~\ref{thm:B}.
In Section~\ref{sec:TRR}, we will discuss some consequences and generalizations
of low genus topological recursion relations. The proofs of the above theorems will be
presented in Section~\ref{sec:proof}.
Our proofs only use topological recursion relations.
This indicates that all these conjectures should be some kind of combinations of topological
recursion relations for all genera.

The author would like to thank H. Xu for presenting a more detailed formulation of Conjecture C.

\section{Equivalence of Conjecture A and B}
\label{sec:A=B}

In \cite{FP}, Faber and Pandharipande considered following Hodge integrals
 over moduli stacks of maps to $V$: For $\gamma_{\alpha_{i}} \in H^{*}(V; \mathbb{C})$,
non-negative integers $k_{i}$ and $b_{j}$,
 \[ \gwih{g, A}{\prod_{i=1}^{n} \tau_{k_{i}}(\gamma_{\alpha_{i}})
    \prod_{j=1}^{m} ch_{b_{j}}(\mathbb{E})}
    := \int_{[\overline{\cal M}_{g,n}(V, A)]^{vir}}
            \prod_{i=1}^{n} \psi_{i}^{k_{i}} \cup ev_{i}^{*}(\gamma_{\alpha_{i}})
                    \cup \prod_{j=1}^{m} ch_{b_{j}}(\mathbb{E}) \]
where $[\overline{\cal M}_{g,n}(V, A)]^{vir}$ is the virtual fundamental cycle in the
moduli space of stable maps from n-marked genus-$g$ stable curves into $V$ with degree
$A \in H^{2}(V; \mathbb{Z})$ (cf. \cite{LiT} and \cite{BF}),
$\psi_{i}$ is the first Chern class of the tautological
line bundle over $\overline{\cal M}_{g,n}(V, A)$ defined by cotangent lines at the $i$-th
marked points on the domain curves, $ev_{i}$ is the evaluation map from the moduli space
to $V$ defined by the image of the $i$-th marked point under stable maps,
$\mathbb{E}$ is the Hodge bundle over the moduli spaces,
and $ch(\mathbb{E})$ is the Chern character of $\mathbb{E}$.
If classes $ch_{b_{j}}(\mathbb{E})$
are omitted in the above expression, this is precisely the descendant Gromov-Witten invariants.
Let $F_{g, \mathbb{E}}(t, s)$ be the generating function for the above
Hodge integrals, where
$t=(t_{n}^{\alpha})$ is the variable for usual Gromov-Witten invariants and
$s=(s_{n})$ is the variable for $ch_{n}(\mathbb{E})$. More precisely
\[ F_{g, \mathbb{E}} := \sum_{A \in H^{2}(V; \mathbb{Z})} q^{A}
            \sum_{m, n \geq 0} \frac{1}{m! \, n!}
            \sum_{\begin{array}{c} \alpha_{1}, \, \cdots,  \, \alpha_{n} \\
                     k_{1}, \cdots, k_{n} \\
                     b_{1}, \cdots, b_{m} \end{array}}
            \prod_{i=1}^{n} t_{k_{i}}^{\alpha_{i}} \prod_{j=1}^{m} s_{b_{j}}
                \gwih{g, A}{\prod_{i=1}^{n} \tau_{k_{i}}(\gamma_{\alpha_{i}})
                     \prod_{j=1}^{m} ch_{b_{j}}(\mathbb{E})}. \]
In particular $F_{g} \, := \, F_{g, \mathbb{E}}\mid_{s=0}  $ is the usual generating function
of genus-$g$ Gromov-Witten invariants. Moreover, define
\[ Z_{\mathbb{E}} := \exp \{ \sum_{g \geq 0} \hbar^{g-1} F_{g, \mathbb{E}} \}\]
where $\hbar$ is a parameter used to separate information from different genera.

For convenience, we identify $\grava{n}$ with
$\frac{\partial}{\partial t_{n}^{\alpha}}$ and $ch_{n}(\mathbb{E})$ with
$\frac{\partial}{\partial s_{n}}$. Moreover we
use the convention that $\grava{0} = \ga$ and $\grava{n}=0$ is $n<0$.
Define
\begin{equation}
D_{2l-1} := - ch_{2l-1}(\mathbb{E})
    - \frac{B_{2l}}{(2l)!} \left\{ \sum_{n, \alpha} \tilde{t}_{n}^{\alpha}
      \,\, \grava{n+2l-1}
    - \frac{\hbar}{2} \sum_{i=0}^{2l-2}
        (-1)^{i} \, \grava{i} \, \gravua{2l-2-i} \right\}
\end{equation}
where $B_{2l}$ are the Bernoulli numbers defined by
\[ \frac{x}{e^{x}-1} = \sum_{r=0}^{\infty} \frac{B_{r}}{r !} x^{r}.
\]

The following formula was proved in \cite{FP}:
\begin{equation} D_{2l-1} Z_{\mathbb{E}} = 0 \end{equation}
for $l \geq 1$.
Taking the coefficients of $\hbar^{g-1}$ in $Z_{\mathbb{E}}^{-1} D_{2l-1} Z_{\mathbb{E}}$
  we  have
\begin{eqnarray}
0 &=& - \frac{(2l)!}{B_{2l}} \gwiigE{ch_{2l-1}(\mathbb{E}) }
    -  \sum_{n, \alpha} \tilde{t}_{n}^{\alpha}
      \gwiigE{ \grava{n+2l-1} }   \nonumber \\
&&  + \frac{1}{2}  \sum_{i=0}^{2l-2}
        (-1)^{i} \left\{ \gwiihE{g-1}{ \grava{i} \gravua{2l-2-i} }
        + \sum_{h = 0}^{g}
        \gwiihE{h}{ \grava{i}} \gwiihE{g-h}{\gravua{2l-2-i} }
    \right\} \nonumber
\end{eqnarray}
for $g \geq 0$ and $l \geq 1$. Here we use $\gwiigE{\cdots}$ to represent derivatives of
$F_{g, \mathbb{E}}$.
Moreover we also adopt the convention that $\gwiihE{g}{\cdots} = 0$ if $g<0$.
This is needed when we consider the genus-0 case of the above formula.
Define
\[
P_{l,\mathbb{E}} := -  \sum_{n, \alpha} \tilde{t}_{n}^{\alpha}
      \grava{n+2l-1}
    + \sum_{i=0}^{2l-2}
        (-1)^{i} \gwiihE{0}{\grava{i}} \gravua{2l-2-i}.
\]
Then the above equation can be rewritten as
\begin{eqnarray}
&& \frac{(2l)!}{B_{2l}} \gwiigE{ch_{2l-1}(\mathbb{E}) }
    - \gwiigE{P_{l,\mathbb{E}}}
    \nonumber \\
&=&  \frac{1}{2}  \sum_{i=0}^{2l-2}
        (-1)^{i} \left\{ \gwiihE{g-1}{ \grava{i} \gravua{2l-2-i} }
        + \sum_{h = 1}^{g-1}
        \gwiihE{h}{ \grava{i}} \gwiihE{g-h}{\gravua{2l-2-i} }
    \right. \nonumber \\
&& \hspace{80pt} \left.
    - \delta_{g, 0} \gwiihE{0}{ \grava{i}} \gwiihE{0}{\gravua{2l-2-i} }
    \right\} \label{eqn:FPE}
\end{eqnarray}
for $g \geq 0$ and $l \geq 1$.

It is also observed in \cite{FP} that
\begin{equation} \label{eqn:vanishE}
 \gwiigE{ch_{2l-1}(\mathbb{E}) } = 0
\end{equation}
if $l > g$. So in this case, if we set $s=0$, equation \eqref{eqn:FPE} implies
 the following relation for
pure Gromov-Witten invariants
\begin{eqnarray}
0
&=&  \gwiig{P_{l}} - \frac{1}{2}  \,\, \delta_{g, 0} \sum_{i=0}^{2l-2}
        (-1)^{i}   \gwii{\grava{i}} \gwii{\gravua{2l-2-i} } \label{eqn:FPGW} \\
&&    + \frac{1}{2}  \sum_{i=0}^{2l-2} (-1)^{i} \left\{ \gwiih{g-1}{ \grava{i} \gravua{2l-2-i} }
         + \sum_{h = 1}^{g-1}
        \gwiih{h}{ \grava{i}} \gwiih{g-h}{\gravua{2l-2-i} }
        \right\} \nonumber
\end{eqnarray}
for $l > g$, where
\begin{equation}
P_{l} := P_{l,\mathbb{E}} \mid_{s =0} \,\, = \, -  \sum_{n, \alpha} \tilde{t}_{n}^{\alpha}
      \grava{n+2l-1}
    + \sum_{i=0}^{2l-2}
        (-1)^{i} \gwii{\grava{i}} \gravua{2l-2-i}.
\end{equation}
For convenience,
we extend the notation $\gwiig{\cdots}$ in such a way that $\gwiig{\vw_{1} \cdots \vw_{k}}$
means the covariant derivative of $F_{g}$ with respect to any vector fields
$\vw_{1}, \cdots , \vw_{k}$ on the big phase space. This notation has been used in the
above equations.
Note that equation \eqref{eqn:FPGW} is just a combination of
conjecture A and conjecture B. It implies that
Conjecture A holds for genus $g \geq 0$ if and only if conjecture B holds for
genus-$(g+1)$. This proves part (a) of Theorem~\ref{thm:B}.

To get more information from equation~\eqref{eqn:FPGW}, we need a better understanding of the vector
field $P_{l}$. For this purpose, we will use the following operator which was studied in \cite{L1}:
\begin{equation} \label{eqn:T}
 T(\vw) := \tau_{+}(\vw) - \gwii{\vw \, \gua} \ga
\end{equation}
for any vector field $\vw$ on the big phase space, where $\tau_{+}(\vw)$ is a linear operator
defined by $\tau_{+}(\grava{n}) := \grava{n+1}$. We will also use
$\tau_{k}(\vw) := \tau_{+}^{k}(\vw)$. The operator $T$ is very useful in the study of
topological recursion relations. For example, a topological recursion relation discovered
by Eguchi-Xiong \cite{EX} can be reformulated as
\begin{equation} \label{eqn:EX}
 \gwiig{ T^{3g-1}(\vw)} = 0
\end{equation}
for any $g>0$ and any vector field $\vw$ (cf. \cite{L1}). This equation follows from
a simple fact that $\psi_{1}^{3k-1}=0$ on the moduli space
$\overline{\cal M}_{g,1}$ due to a dimension count.

Let
\[ \vs:= - \sum_{n, \alpha} \tilde{t}_{n}^{\alpha} \grava{n-1} \]
be the {\it string vector field}.
The following lemma follows from \cite[Equation (33)]{L2}. It also follows
from equation~\eqref{eqn:Ttau2} below and the dilaton equation.
\begin{lem} \label{lem:PS}
For $l \geq 1$,
\[ P_{l} = T^{2l}(\vs). \]
\end{lem}
Therefore combining equations \eqref{eqn:FPGW} and \eqref{eqn:EX}, we have
\[
    \sum_{i=0}^{2l-2} (-1)^{i} \left\{ \gwiih{g-1}{ \grava{i} \gravua{2l-2-i} }
         + \sum_{h = 1}^{g-1}
        \gwiih{h}{ \grava{i}} \gwiih{g-h}{\gravua{2l-2-i} }
        \right\} = 0
\]
if $g>0$ and $l \geq \max\{g+1, \frac{3g-1}{2} \}$.
Replacing $l$ by $k+1$, we have
\begin{lem} \label{lem:FPGWEX}
For $g>0$ and $k \geq \max\{g, \frac{3g-3}{2} \}$,
\[
   \sum_{i=0}^{2k} (-1)^{i} \left\{ \gwiih{g-1}{ \grava{i} \gravua{2k-i} }
         + \sum_{h = 1}^{g-1}
        \gwiih{h}{ \grava{i}} \gwiih{g-h}{\gravua{2k-i} }
        \right\} =0.
\]
\end{lem}
Note that this lemma is a special case of
Theorem~\ref{thm:A} and Theorem~\ref{thm:allg}.
Note that when $g=1$,  Lemma \ref{lem:FPGWEX} is
 equivalent to a genus-0 equation
\begin{equation} \label{eqn:g0conjA}
 \sum_{i=0}^{2k} (-1)^{i} \gwii{ \grava{i} \,\, \gravua{2k-i} } = 0
\end{equation}
for $k>0$. This is precisely the genus-0 case of conjecture A.

\section{Low genus TRR}
\label{sec:TRR}

For simplicity, "topological recursion relations" will be abbreviated as TRR.
Recall the {\it genus-0 TRR}:
\[ \gwii{\tau_{1}(\vw) \, \, \vu \, \vv}
    = \gwii{\vw \, \gua} \gwii{\ga \, \vu \, \vv} \]
for any vector fields $\vw$, $\vu$ and $\vv$.
As observed by Witten \cite{W}, this formula implies the {\it generalized WDVV equation}:
\[ \gwii{\vw_{1} \vw_{2} \gua} \gwii{\ga \, \vw_{3} \vw_{4}}
    = \gwii{\vw_{1} \vw_{3} \gua} \gwii{\ga \, \vw_{2} \vw_{4}}. \]

We now prove some other useful consequences of the genus-0 TRR.
First, we prove the following formula
\begin{lem} \label{lem:genSrec}
For $m \geq 0$,
\begin{eqnarray}
&&  \sum_{i=0}^{m} (-1)^{i} \gwii{\vw \, \gravua{i}} \gwii{\grava{m-i} \,  \vv}
        \nonumber \\
&=& \gwii{ \tau_{m+1}(\vw) \, \, \vv} + (-1)^{m} \gwii{ \vw \, \tau_{m+1}(\vv)}
        \label{eqn:genSrec}
\end{eqnarray}
for any vector fields $\vw$ and $\vv$.
\end{lem}
{\bf Proof}:
We first note that when $m=0$ this formula has the following form
\begin{equation} \label{eqn:Srec}
  \gwii{\vw \, \gua} \gwii{\ga \,  \vv}
    = \gwii{ \tau_{1}(\vw) \, \, \vv} +  \gwii{ \vw \, \tau_{1}(\vv)}.
\end{equation}
This is exactly \cite[Equation (10)]{LT}, which was proved using
the string equation and the genus-0 TRR.
When $m=1$, this formula is exactly \cite[Lemma 4.3 (iii)]{LT}.
We now prove the lemma by induction on $m$. Assume the lemma is true for $m=k$.

For $m=k+1$, we apply equation \eqref{eqn:Srec} to obtain the following formula
for the left hand side of equation \eqref{eqn:genSrec}:
\begin{eqnarray*}
{\rm LHS} &=& \sum_{i=0}^{k} (-1)^{i} \gwii{\vw \, \gravua{i}}
            \left\{ - \gwii{\grava{k-i} \, \tau_{1}(\vv)}
                         + \gwii{\grava{k-i} \, \gub} \gwii{\gb \,  \vv} \right\} \\
&&         + (-1)^{k+1} \gwii{ \vw \, \grava{k+1}} \gwii{\ga \,  \vv}
\end{eqnarray*}
Applying the induction hypothesis to the first two terms, we obtain
\begin{eqnarray*}
{\rm LHS} &=& - \left\{ \gwii{\tau_{k+1}(\vw) \, \tau_{1}(V)} + (-1)^{k} \gwii{\vw \, \tau_{k+2}(\vv)}
            \right\}
            + \gwii{\tau_{k+1}(\vw) \, \gub} \gwii{\gb \,  \vv}.
\end{eqnarray*}
The lemma follows by applying equation \eqref{eqn:Srec} to the last term.
$\Box$

An immediate consequence of Lemma \ref{lem:genSrec} is the following
\begin{lem}[Generalized genus-0 TRR] \label{lem:geng0TRR} For $m \geq 0$,
\[ \sum_{i=0}^{m} (-1)^{i} \gwii{\vw \, \gravua{i}} \gwii{\grava{m-i} \, \vu \, \vv}
    = \gwii{ \tau_{m+1}(\vw) \, \,\vu \, \vv} \]
for any vector fields $\vw$, $\vu$ and $\vv$.
\end{lem}
{\bf Proof}:
When $m=0$, this formula is precisely the genus-0 TRR.
Now assume $m > 0$.
Applying the genus-0 TRR to the left hand side of this formula, we obtain
\begin{eqnarray*}
{\rm LHS} &=& \sum_{i=0}^{m-1} (-1)^{i} \gwii{\vw \, \gravua{i}}
                    \gwii{\grava{m-i-1} \, \gub} \gwii{\gb \, \vu \, \vv} \\
    && + (-1)^{m} \gwii{\vw \, \gravua{m}} \gwii{\ga \, \vu \, \vv}
\end{eqnarray*}
Applying Lemma \ref{lem:genSrec} to the first term, we have
\[
{\rm LHS} =  \gwii{\tau_{m}(\vw) \, \gub}  \gwii{\gb \, \vu \, \vv} =
    \gwii{\tau_{m+1}(\vw) \,  \vu \, \vv}
\]
by the genus-0 TRR. The lemma is proved. $\Box$

Note that a special case of Theorem~\ref{thm:allg} is the following:
If $k, n \geq 2$ and $m \geq k+n-3$,
\begin{eqnarray}
 \sum_{i=0}^{m} (-1)^{i} \gwii{\vw_{1} \, \cdots \, \vw_{n} \, \, \gravua{i}}
        \gwii{  \grava{m-i} \, \vv_{1}\, \cdots \, \vv_{k}}
    &=& 0 \label{eqn:conjCg0pt3}
\end{eqnarray}
for any vector fields $\vw_{1}, \ldots, \vw_{n}$ and $\vv_{1}, \ldots, \vv_{k}$.
Therefore if we take derivatives of the formula in Lemma~\ref{lem:geng0TRR}, we obtain
\begin{equation} \label{eqn:conjCg0pt2}
 \sum_{i=0}^{m} (-1)^{i} \gwii{\vw \, \gravua{i}} \gwii{\grava{m-i} \, \vv_{1} \, \cdots \vv_{k}}
    = \gwii{ \tau_{m+1}(\vw) \, \,\vv_{1} \, \cdots \vv_{k}}
\end{equation}
for $k \geq 2$, $m \geq k-2$, and any vector fields $\vw, \,\,\, \vv_{1}, \ldots, \vv_{k}$.

Note that together with genus-0 part of Conjecture B and Lemma~\ref{lem:genSrec},
equations \eqref{eqn:conjCg0pt3} and \eqref{eqn:conjCg0pt2} covers most cases of genus-0
part of Conjecture C.
The remaining cases follows from derivatives of Conjecture B.

We now study analogous results for genus-1 Gromov-Witten invariants.
The {\it genus-1 TRR} is the following:
\[ \gwiione{\tau_{1}(\vw)} = \gwii{\vw \, \gua} \gwiione{\ga}
    + \frac{1}{24} \gwii{ \vw \, \gua  \, \ga} \]
for any vector field $\vw$.
Similar to Lemma \ref{lem:geng0TRR}, we also have the following
\begin{lem}[Generalized genus-1 TRR] \label{lem:geng1TRR} For $m \geq 0$,
\[ \sum_{i=0}^{m} (-1)^{i} \gwii{\vw \, \gravua{i}} \gwiih{1}{\grava{m-i}}
    = \gwiih{1}{ \tau_{m+1}(\vw)} -  \delta_{m, 0} \,\, \frac{1}{24} \gwii{\vw \, \,\gua \, \ga} \]
for any vector field $\, \vw$.
\end{lem}
{\bf Proof}:
When $m=0$, this formula is precisely the genus-1 TRR.
Now assume $m > 0$.
Applying the genus-1 TRR to the left hand side of this formula, we obtain
\begin{eqnarray*}
{\rm LHS} &=& \sum_{i=0}^{m-1} (-1)^{i} \gwii{\vw \, \gravua{i}}
                    \left\{ \gwii{ \tau_{m-1-i}(\ga) \, \gub} \gwiih{1}{\gb} +
                            \frac{1}{24} \gwii{\tau_{m-1-i}(\ga) \, \,\gub \, \gb} \right\} \\
    && + (-1)^{m} \gwii{\vw \, \gravua{m}} \gwiih{1}{\ga }
\end{eqnarray*}
Applying Lemma \ref{lem:genSrec} to the first term and Lemma \ref{lem:geng0TRR}
to the second term, we have
\[
{\rm LHS} =  \gwii{\tau_{m}(\vw) \, \gub} \gwiih{1}{\gb} +
                            \frac{1}{24} \gwii{\tau_{m}(\vw) \, \,\gub \, \gb}
    =  \gwiih{1}{\tau_{m+1}(\vw)}
\]
by the genus-1 TRR. This proves the lemma. $\Box$

Note that a special case of Theorem~\ref{thm:allg} is the following:
For $n \geq 2$, $k \geq 0$ and $m \geq k+n$,
\begin{eqnarray}
 \sum_{i=0}^{m} (-1)^{i} \gwii{\vw_{1} \, \cdots \, \vw_{n} \, \, \gravua{i}}
        \gwiione{  \grava{m-i} \, \vv_{1}\, \cdots \, \vv_{k}}
    &=& 0 \label{eqn:conjCg1pt3}
\end{eqnarray}
for any vector fields $\vw_{1}, \ldots, \vw_{n}$ and $\vv_{1}, \ldots, \vv_{k}$.
Therefore if we take derivatives of the formula in Lemma~\ref{lem:geng1TRR}, we obtain
\begin{equation} \label{eqn:conjCg1pt2}
 \sum_{i=0}^{m} (-1)^{i} \gwii{\vw \, \gravua{i}} \gwiione{\grava{m-i} \, \vv_{1} \, \cdots \vv_{k}}
    = \gwiione{ \tau_{m+1}(\vw) \, \,\vv_{1} \, \cdots \vv_{k}}
\end{equation}
for $k \geq 0$, $m \geq k+1$, and any vector fields $\vw, \,\,\, \vv_{1}, \ldots, \vv_{k}$.

Note that equations~\eqref{eqn:conjCg1pt3} and \eqref{eqn:conjCg1pt2} are
 stronger than corresponding cases of genus-1 part of
Conjecture C in the sense that we do not sum over genus. It implies most cases of the genus-1 part
of conjecture C. The remaining cases follow from genus-1 part of Conjecture B and
its derivatives.

Using the operator $T$, we can reformulate the genus-0 TRR as
\[ \gwii{T(\vw_{1}) \, \vw_{2} \, \vw_{3}} = 0 \]
and the genus-1 TRR as
\[ \gwiione{T(\vw)} = \frac{1}{24} \gwii{\vw \, \gua \, \ga} \]
for any vector fields $\vw$ and $\vw_{i}$.
Replacing $\vw$ by $T^{k}(\vw)$ in the genus-1 TRR, we obtain
\[ \gwiione{T^{k}(\vw)} = 0 \]
for all vector field $\vw$ if $k \geq 2$.

For later use, we also recall the genus-2 Mumford relation (cf. \cite{Ge}) as
formulated in \cite{L1}:
\begin{eqnarray}
\gwiitwo{T^{2}(\vw)} &=&
    \frac{7}{10} \gwiione{\ga} \gwiione{\{\gua \circ \vw\}}
    + \frac{1}{10} \gwiione{\ga \, \{\gua \circ \vw\}}  \nonumber \\
&&     - \frac{1}{240} \gwiione{\vw \, \{\ga \circ \gua\}}
    + \frac{13}{240} \gwii{\vw \, \ga \, \gua \, \gub}
        \gwiione{\gb} \nonumber \\
&&    + \frac{1}{960} \gwii{ \vw \, \gua \, \ga \, \gub \, \gb}
    \label{eqn:TRRg2M}
\end{eqnarray}
for any vector field $\vw$. Here we have used the quantum product for vector fields
on the big phase space defined by
\[ \vw_{1} \circ \vw_{2} := \gwii{\vw_{1} \, \vw_{2} \, \gua} \, \ga. \]
Basic properties of this product can be found in \cite{L1}.
Replacing $\vw$ by $T^{i}(\vw)$ in equation \eqref{eqn:TRRg2M} and using genus-0 and genus-1
TRR, we obtain
\begin{eqnarray}
 \gwiitwo{ T^{3}(\vw) }
&=& \frac{1}{20} \gwiione{\left\{ \vw \qp \gua \qp \ga \right\} }
    + \frac{1}{1152} \gwii{ \vw \,\, \gua \,\, \ga \,\, \left\{ \gub \qp \gb \right\}} \nonumber \\
&&    + \frac{1}{480} \gwii{\left\{ \vw \qp \gua \right\} \, \, \ga \,\, \gub \,\, \gb}
    \label{eqn:g2T3},
\end{eqnarray}
\begin{eqnarray}
 \gwiitwo{ T^{4}(\vw) }
 &=& \frac{1}{1152} \gwii{\left\{ \vw \qp \gua \qp \ga \right\} \,\, \gub \, \, \gb},
    \hspace{120pt}
    \label{eqn:g2T4}
\end{eqnarray}
and
\[ \gwiitwo{T^{k}(\vw)} = 0 \hspace{290pt} \]
for any vector field $\vw$ if $k \geq 5$.

To prove Theorem~\ref{thm:B}, we also need to use the following genus-3 equation (cf. \cite{KL}):
{\allowdisplaybreaks
\begin{eqnarray}
&& \gwiih{3}{T^{3}(\vw)} \nonumber \\
&=& - \frac{1}{252} \gwiitwo{ \vw \, T(\ga \qp \gua)}
    + \frac{5}{42} \gwiitwo{T(\ga) \, \{ \vw \qp \gua \}} \nonumber \\
&&  + \frac{13}{168} \gwiitwo{ T(\gua) } \gwii{ \ga \, \vw \, \gub \, \gb}
    + \frac{41}{21} \gwiitwo{ T(\gua) } \gwiione{ \left\{\ga \qp \vw\right\}}
        \nonumber \\
&&  - \frac{13}{168} \gwiitwo{ \left\{\vw \qp \ga \qp \gua\right\} }
    + \frac{1}{280} \gwiione{ \vw \gua} \gwiione{ \ga \, \left\{\gub \qp \gb\right\}}
        \nonumber \\
&&  - \frac{23}{5040} \gwiione{\gua} \gwiione{\ga \vw \left\{\gub \qp \gb\right\}}
 - \frac{47}{5040} \gwiione{\gua} \gwiione{\ga \gub} \gwii{\gb \, \vw \, \gum \, \gm}
        \nonumber \\
&&  - \frac{5}{1008} \gwiione{\vw \, \gua} \gwii{\ga \gub \gb \gum} \gwiione{\gm}
  + \frac{23}{504} \gwiione{\gua} \gwii{\ga \vw \gub \gb \gum} \gwiione{\gm}
         \nonumber \\
&&   + \frac{11}{140} \gwiione{\gua \gub} \gwiione{\ga \, \left\{\gb \qp \vw\right\}}
  - \frac{4}{35} \gwiione{\gua} \gwiione{\ga \gub} \gwiione{\left\{\gb \qp \vw\right\}}
            \nonumber \\
&&  + \frac{2}{105} \gwiione{\vw \, \gua} \gwiione{\left\{\ga \qp \gb\right\}} \gwiione{\gub}
  + \frac{89}{210} \gwiione{\gua} \gwii{ \ga \, \vw \, \gub \, \gum} \gwiione{\gb} \gwiione{\gm}
            \nonumber\\
&& - \frac{1}{210} \gwiione{\gua} \gwiione{\ga \gub \left\{\gb \qp \vw\right\}}
    + \frac{1}{140} \gwiione{\vw \, \gua \gub} \gwiione{\left\{\ga \qp \gb\right\}}
            \nonumber \\
&& + \frac{23}{140} \gwiione{\gua \gub} \gwii{\ga \gb \vw \gum} \gwiione{\gm}
    - \frac{3}{140} \gwiione{\gua \gub} \gwiione{\left\{\ga \qp \gb\right\} \vw}
            \nonumber \\
&& - \frac{1}{4480} \gwiione{\vw \gua} \gwii{\ga \gb \gub \gm \gum}
    + \frac{13}{8064} \gwiione{\gua} \gwii{\ga \vw \gub \gb \gum \gm}
    \nonumber \\
&& - \frac{1}{2240} \gwiione{\vw \gua \gub} \gwii{\ga \gb \gum \gm}
    + \frac{41}{6720} \gwiione{\gua \gub} \gwii{\ga \gb \vw \gum \gm}
    \nonumber \\
&& + \frac{1}{53760} \gwii{\vw \gua \ga \gub \gb \gum \gm}
    - \frac{1}{210} \gwiione{\left\{\vw \qp \gua\right\}} \gwiione{\ga \gub \gb}
    \nonumber \\
&& - \frac{1}{5760} \gwiione{\vw \gua \ga \left\{\gub \qp \gb\right\} }
    - \frac{1}{2688} \gwiione{\gua \ga \gub} \gwii{\gb \vw \gum \gm}
    \nonumber \\
&& - \frac{1}{5040} \gwiione{\gua \ga \gub \left\{\gb \qp \vw\right\}}
   + \frac{1}{3780} \gwiione{ \vw \ga \gb \gm } \gwii{\gua \gub \gum}
    \nonumber \\
&&    + \frac{1}{252} \gwiione{ \ga \gb \gm } \gwii{ \vw \gua \gub \gum}.
    \label{eqn:g3TRR}
\end{eqnarray}}
Replacing $\vw$ by $T^{3}(\vw)$ in this equation, we get
\begin{eqnarray}
\gwiih{3}{T^{6}(\vw)}
&=& \frac{7}{5760} \gwiione{\left\{\vw \qp \Delta \qp \Delta\right\}}
    + \frac{11}{2903040} \gwii{\vw \, \Delta \, \Delta \, \Delta }
        \nonumber \\
&& + \frac{19}{967680} \gwii{\left\{\vw \qp \Delta \right\} \, \Delta \, \gua \, \ga}
        + \frac{1}{120960} \gwii{\vw  \,
            \left\{\Delta \qp \Delta \right\} \, \gua \, \ga}
     \nonumber \\
&&
    + \frac{1}{60480} \gwii{\left\{\vw \qp \gua \right\} \, \ga \, \Delta \, \Delta}
        \nonumber \\
&&
+ \frac{1}{11520} \gwii{\left\{\vw \qp \Delta \qp \gua \right\} \, \ga \, \gub \, \gb}
    \label{eqn:g3T6}
\end{eqnarray}
for any vector field $\vw$.

\section{Proof of main theorems}
\label{sec:proof}

\allowdisplaybreaks

In order to see the connection between TRR and conjectures in the introduction, we give a
new formulation to these conjectures and prove that these conjectures are correct for all genera if
$m$ is sufficiently large. The key point here is that we will use the operator $T$ instead of
descendant operator $\tau$. The following lemma gives the relation between $T^{k}(\vw)$
and $\tau_{k}(\vw)$.

\begin{lem} \label{lem:Ttau}
For any $k \geq 0$ and any vector field $\vw$ on the big phase space,
\[ \tau_{k}(\vw) = T^{k}(\vw) + \sum_{i=0}^{k-1} \gwii{\tau_{k-1-i}(\vw) \, \gua} T^{i}(\ga). \]
\end{lem}
{\bf Proof}: We prove the lemma by induction on $k$.
If $k=0$, it is trivial. If $k=1$, it is just the definition of $T$.
Assume it is true for $k=r$, then for $k=r+1$, we have
\[
\tau_{k}(\vw) = \tau_{+}(\tau_{k-1}(\vw))
    = T(\tau_{k-1}(\vw)) + \gwii{\tau_{k-1}(\vw) \, \gua} \ga
\]
by the definition of $T$. By the induction hypothesis,
\begin{eqnarray*}
\tau_{k}(\vw) &=& T \left(
        T^{k-1}(\vw) + \sum_{i=0}^{k-2} \gwii{\tau_{k-2-i}(\vw) \, \gua} T^{i}(\ga) \right)
            + \gwii{\tau_{k-1}(\vw) \, \gua} \ga \\
&=& T^{k}(\vw) + \sum_{i=0}^{k-1} \gwii{\tau_{k-1-i}(\vw) \, \gua} T^{i}(\ga).
\end{eqnarray*}
The lemma is proved.
$\Box$

The following proposition shows that in certain combinations,  $\tau_{k}$ can be replaced by
$T^{k}$.
\begin{pro} \label{pro:Ttau}
For any contravariant tensors $P$ and $Q$ on the big phase space and $m \geq 0$,
\[ \sum_{j=0}^{m} (-1)^{j} P(\gravua{j}) \,\,  Q(\grava{m-j})
    = \sum_{j=0}^{m} (-1)^{j} P(T^{j}(\gua)) \, \, Q(T^{m-j}(\ga)).\]
\end{pro}
{\bf Proof}:
By Lemma \ref{lem:Ttau}, we have
\[ P(\gravua{j})= P(T^{j}(\gua)) + \sum_{r=0}^{j-1} \gwii{\tau_{j-1-r}(\gua) \, \gub} P(T^{r}(\gb)).\]
We also have a similar formula for $Q(\grava{m-j})$. Therefore
the difference of the two sides of the equation in the proposition
is given by
\begin{eqnarray*}
&& {\rm LHS} - {\rm RHS} \\
&=& \sum_{j=0}^{m} (-1)^{j} \sum_{s=0}^{m-j-1} \gwii{\grava{m-j-1-s} \gum}
                P(T^{j}(\gua)) Q(T^{s}(\gm))  \\
&& + \sum_{j=0}^{m} (-1)^{j} \sum_{r=0}^{j-1} \gwii{\gravua{j-1-r} \gub}
                P(T^{r}(\gb)) Q(T^{m-j}(\ga)) \\
&& + \sum_{j=0}^{m} (-1)^{j} \sum_{r=0}^{j-1} \sum_{s=0}^{m-j-1}
            \gwii{\gravua{j-1-r} \gub} \gwii{\grava{m-j-1-s} \gum}
                P(T^{r}(\gb)) Q(T^{s}(\gm)).
\end{eqnarray*}
The third term can be written as
\begin{eqnarray*}
&&   \sum_{r=0}^{m-1} \sum_{s=0}^{m-1} P(T^{r}(\gb)) Q(T^{s}(\gm))
        \sum_{j=r+1}^{m-1-s} (-1)^{j}
            \gwii{\gravua{j-1-r} \gub} \gwii{\grava{m-j-1-s} \gum} \\
&=& \sum_{r=0}^{m-1} \sum_{s=0}^{m-1} P(T^{r}(\gb)) Q(T^{s}(\gm)) \\
&& \hspace{60pt} \cdot
        (-1)^{r+1} \left\{
            \gwii{\gravub{m-r-s-1} \, \gum}
                + (-1)^{m-r-s} \gwii{\gub \, \gravum{m-r-s-1}} \right\}
\end{eqnarray*}
by Lemma \ref{lem:genSrec}.
Hence it can be canceled with the other two terms in the above expression.
The proposition is thus proved.
$\Box$

As a consequence of this proposition, we have the following
\begin{thm} \label{thm:reform}
In conjectures A, B, C in the introduction,
$\gravua{j}$ and $\grava{m-j}$ can be replaced by $T^{j}(\gua)$ and $T^{m-j}(\ga)$
respectively.
\end{thm}

We are now ready to prove main theorems of this paper.

\vspace{10pt}
\noindent
{\bf Proof of Theorem \ref{thm:allg}}:

By Proposition~\ref{pro:Ttau},
\begin{eqnarray}
&& \sum_{j=0}^{m}
      (-1)^{j} \gwiih{h}{\grava{j}  \, \, \vw_{1} \, \cdots \, \vw_{a}}
    \gwiih{g-h}{\gravua{m-j} \,\, \vv_{1} \, \cdots \, \vv_{b} } \nonumber \\
&=& \sum_{j=0}^{m}
      (-1)^{j} \gwiih{h}{T^{j}(\ga) \,\, \vw_{1} \, \cdots \, \vw_{a}}
    \gwiih{g-h}{T^{m-j}(\gua) \,\, \vv_{1} \, \cdots \, \vv_{b}}
    \label{eqn:Ttauhg}
\end{eqnarray}
for any vector fields $\vw_{i}$ and $\vv_{j}$.
Since the dimension of moduli space $\overline{\cal M}_{g, k}$ is $3g-3+k$,
$\psi_{1}^{j} = 0$ on $\overline{\cal M}_{g, k}$ if $j>3g-3+k$.
Translating this fact to a universal equation for Gromov-Witten invariants, we get
\[ \gwiig{T^{j}(\vv) \, \vw_{1} \, \cdots \, \vw_{k}} = 0 \]
 for any vector fields $\vv$ and $\vw_{i}$ if $j>3g-2+k$.
 In the genus-0 case, we also require $k \geq 2$ since $\overline{\cal M}_{g, n}$
 does not exist if $n<3$.
 So any term in equation \eqref{eqn:Ttauhg} which is non-zero requires that
$j \leq 3h-2+a$ and $m-j \leq 3(g-h)-2+b$. But this requires that $m \leq 3g+a+b-4$.
Equation~\eqref{eqn:allg1} in Theorem~\ref{thm:allg} is thus proved.

To prove equation~\eqref{eqn:allg2}, we first prove the special case that $b=0$.
Note that equation~\eqref{eqn:Srec} can be written in the following form
\begin{equation} \label{eqn:Ttaug0}
\gwii{T(\vw) \,\, \vv} = - \gwii{\vw \,\, \tau_{+}(\vv)}
\end{equation}
for all vector fields $\vw$ and $\vv$.
So by Proposition~\ref{pro:Ttau}, Lemma~\ref{lem:Ttau} can be rewritten as
\begin{equation} \label{eqn:Ttau2}
T^{k}(\vw) = \tau_{k}(\vw) - \sum_{i=0}^{k-1} (-1)^{i} \gwii{\vw \, \tau_{i}(\gua)} \tau_{k-1-i}(\ga).
\end{equation}
Note that this is exactly \cite[Equation (32)]{L2}. We also notice that if
we replace $\vw$ in this equation by the dilaton vector field, then we obtain
Lemma~\ref{lem:PS} by the dilaton equation.

By equation~\eqref{eqn:EX},
$\gwiig{T^{k}(\vw)} = 0$ if $k \geq 3g-1$. So
by equation~\eqref{eqn:Ttau2} with $m = k-1$, we have
\begin{equation} \label{eqn:genTRR}
\sum_{i=0}^{m} (-1)^{i} \gwii{\vw \, \gravua{i}} \gwiig{\grava{m-i}}
    = \gwiig{\tau_{m+1}(\vw)}
\end{equation}
for any vector field $\vw$ if $g \geq 1$ and $m \geq 3g-2$.
This is a special case of equation~\eqref{eqn:allg2} with $b=0$.
We can view this equation as a generalization of Lemma \ref{lem:geng1TRR}.
Taking derivatives of this equation and using equation~\eqref{eqn:allg1},
we obtain equation~\eqref{eqn:allg2}. Theorem~\ref{thm:allg} is thus proved.
$\Box$

\vspace{10pt}
\noindent
{\bf Proof of Theorem \ref{thm:A}}:

The genus-0 part of Theorem \ref{thm:A} is precisely equation \eqref{eqn:g0conjA}.
To illustrate the idea that Conjecture A should follow from topological recursion relations,
we provide another proof here. Since Conjecture A is trivial when $m$ is odd,
we assume that $m$ is an even positive integer.
By an argument similar to the proof of Proposition~\ref{pro:Ttau} and the genus-0 TRR, we have
\begin{eqnarray}
\sum_{j=0}^{m} (-1)^{j} \gwii{\grava{j} \,\, \vw \,\, \gravua{m-j}}
&=& \sum_{j=0}^{m} (-1)^{j} \gwii{T^{j}(\ga) \,\, \vw \,\, T^{m-j}(\gua)} \,\, = \,\, 0
\label{eqn:derg0conjA}
\end{eqnarray}
for any vector field $\vw$.

Now we consider the special case that $\vw$ is the {\it Euler vector field}
\[ \vx := - \sum_{n, \alpha} \left(n+b_{\alpha}  - \frac{3-d}{2}\right)
    \tilde{t}_{n}^{\alpha} \grava{n} - \sum_{n, \alpha, \beta} {\cal C}_{\alpha}^{\beta}
            \tilde{t}_{n}^{\alpha} \gravb{n-1} \]
where $d= \frac{1}{2} \, {\rm dim}_{\mathbb{R}}(V)$, $b_{\alpha}= \frac{1}{2} \left\{ {\rm deg}(\ga) -d +1 \right\}$,
and ${\cal C}_{\alpha}^{\beta}$ represents the multiplication
by the first Chern class of $V$, i.e. $c_{1}(V) \cup \ga = {\cal C}_{\alpha}^{\beta} \gb$.
By \cite[Lemma 1.4 (3)]{LT}, which follows from the quasi-homogeneity equation, we have
\begin{eqnarray*}
 \sum_{j=0}^{m} (-1)^{j} \gwii{\grava{j} \,\, \vx \,\, \gravua{m-j}}
&=& \sum_{j=0}^{m} (-1)^{j} \left\{ (m+1) \gwii{\grava{j} \,\, \gravua{m-j}} \right. \\
&& \hspace{60pt} \left.
            +  {\cal C}_{\alpha}^{\mu} \eta_{\mu \nu} \gwii{\gravun{j-1} \, \gravua{m-j}} \right. \\
&& \hspace{60pt} \left.
            +  \eta^{\alpha \beta} {\cal C}_{\beta}^{\mu}
                    \gwii{\grava{j} \, \gravm{m-j-1}}
            \right\}.
\end{eqnarray*}
Since both ${\cal C}_{\alpha}^{\mu} \eta_{\mu \nu}$ and
$\eta^{\alpha \beta} {\cal C}_{\beta}^{\mu}$ are (super)symmetric with respect to the free indices,
the last two terms are equal to 0. Therefore we have
\[  \sum_{j=0}^{m} (-1)^{j} \gwii{\grava{j} \,\, \vx \,\, \gravua{m-j}}
=  (m+1) \sum_{j=0}^{m} (-1)^{j}   \gwii{\grava{j} \,\, \gravua{m-j}}. \]
So the genus-0 part of Theorem \ref{thm:A} follows from equation
\eqref{eqn:derg0conjA} with $\vw$ replaced by $\vx$.

For $g>0$, we have
\[ \gwiig{T^{j}(\vw) \, T^{m-j}(\vv)} = 0 \]
for any vector fields $\vw$ and $\vv$ if $m \geq 3g$. This follows from the fact
that $\psi_{1}^{j} \psi_{2}^{m-j} = 0$ on $\overline{\cal M}_{g, 2}$ since
the dimension of $\overline{\cal M}_{g, 2}$ is $3g-1$.
Therefore Theorem~\ref{thm:A} follows from Proposition~\ref{pro:Ttau}.
$\Box$.

\vspace{10pt}
\noindent
{\bf Proof of Theorem \ref{thm:B}}:

Part (a) of Theorem \ref{thm:B} was proved in Section~\ref{sec:A=B}.
 When $g=k=0$, conjecture B is just the genus-0 string equation.
 Moreover, when $g=0$ and $l>0$, equation \eqref{eqn:FPGW} is exactly the genus-0 part of
conjecture B with $k=l$.
By Lemma~\ref{lem:PS}, for genus bigger than 0,
conjecture B can be reformulated as
\begin{equation} \label{eqn:Bphi}
\gwiig{T^{2k}(\vs)} + \frac{1}{2} \Phi_{g, k-1} = 0
\end{equation}
for $k \geq g$, where $\vs$ is the string vector field and
\[ \Phi_{g, k} := \sum_{j=0}^{2k} (-1)^{j} \sum_{h=1}^{g-1} \gwiih{h}{\gravua{j}}
        \gwiih{g-h}{\grava{2k-j}}.\]
Since $\gwiione{T^{2}(\vw)}=0$  for any vector field $\vw$ by the genus-1 TRR,
genus-1 case of Conjecture B follows from equation~\eqref{eqn:Bphi} rather trivially.
So we only need to consider the case of $g >1$.
Part (b) of Theorem~\ref{thm:B} follows from equation~\eqref{eqn:EX} and
Theorem~\ref{thm:allg}. Part (c) of Theorem~\ref{thm:B} follows from
part (b) except the cases $k=g=2$ and $k=g=3$.

We first look at the genus-2 case. By Proposition~\ref{pro:Ttau}
and genus-1 TRR,
\begin{eqnarray*}
\Phi_{2, 1} &=& \sum_{j=0}^{2} (-1)^{j}  \gwiih{1}{T^{j}(\gua)} \gwiih{1}{T^{2-j}(\ga)}
\,\, =\,\,  -  \gwiih{1}{T(\gua)} \gwiih{1}{T(\ga)}\\
&=& - \frac{1}{576} \gwii{\gua \, \gub \, \gb}
                \gwii{\ga \, \gum \, \gm}
\,\, = \,\, - \frac{1}{576} \gwii{\Delta \, \gub \, \gb}
\end{eqnarray*}
where \[ \Delta := \gua \circ \ga.\]
On the other hand, since $\vs \circ \Delta = \Delta$,
by equation \eqref{eqn:g2T4},
\[ \gwiitwo{T^{4}(\vs)} = \frac{1}{1152} \gwii{ \Delta \, \gua \, \ga}. \]
By equation \eqref{eqn:Bphi},
this proves the $g=k=2$ case of conjecture B.

Now we consider the $g=k=3$ case.
By Proposition~\ref{pro:Ttau} and the genus-1 TRR, we have
\begin{eqnarray*}
\frac{1}{2} \Phi_{3, 2}
&=& \sum_{j=0}^{4} (-1)^{j}  \gwiih{1}{T^{j}(\gua)} \gwiih{2}{T^{4-j}(\ga)} \\
&=& \gwiih{1}{\gua} \gwiih{2}{T^{4}(\ga)} - \gwiih{1}{T(\gua)} \gwiih{2}{T^{3}(\ga)} .
\end{eqnarray*}
By the genus-1 TRR and equations \eqref{eqn:g2T3} and \eqref{eqn:g2T4},
\begin{eqnarray*}
\frac{1}{2} \Phi_{3, 2}
&=& \frac{1}{1152} \gwiione{\gua} \gwii{(\ga \circ \Delta) \, \gub \, \gb}  \\
&& - \frac{1}{24} \gwii{\gua \, \gum \, \gm}
        \left\{ \frac{1}{20} \gwiione{(\ga \circ \Delta)}
            + \frac{1}{1152} \gwii{\ga \, \gub \, \gb \, \Delta} \right. \\
&& \hspace{100pt}
       \left.     + \frac{1}{480} \gwii{(\ga \circ \gub) \, \gb \, \gun \, \gn}
                \right\} \\
&=& - \frac{7}{5760} \gwiione{(\Delta \circ \Delta)}
    - \frac{1}{27648} \gwii{\Delta \, \Delta \, \gua \, \ga}
    - \frac{1}{11520} \gwii{(\Delta \circ \gua) \, \ga \, \gub \, \gb}.
\end{eqnarray*}
On the other hand, since $\vs \circ \ga = \ga$ and
$\gwii{\vs \, \ga \, \gb\, \gm}=0$ for all $\alpha$, $\beta$ and $\mu$,
by equation~\eqref{eqn:g3T6},  we have
\begin{eqnarray*}
\gwiih{3}{T^{6}(\vs)} &=&  \frac{7}{5760} \gwiione{(\Delta \circ \Delta)}
    + \frac{1}{27648} \gwii{\Delta \, \Delta \, \gua \, \ga}
    + \frac{1}{11520} \gwii{(\Delta \circ \gua) \, \ga \, \gub \, \gb}.
\end{eqnarray*}
 So the $g=k=3$ case of Conjecture B follows from
equation \eqref{eqn:Bphi}.
This finishes the proof of Theorem~\ref{thm:B}.
$\Box$


\vspace{30pt} \noindent
Department of Mathematics  \\
University of Notre Dame \\
Notre Dame,  IN  46556, USA \\

\vspace{10pt} \noindent E-mail address: {\it xliu3@nd.edu}

\end{document}